\documentclass[12pt]{article}

\usepackage{epsf}
\usepackage{psfrag}
\usepackage{amsthm}
\usepackage{amsmath}

\renewcommand{\epsfsize}[2]{0.6\textwidth}

\theoremstyle{plain}
\newtheorem{thm}{Theorem}[section]
\newtheorem{cor}[thm]{Corollary}

\newtheorem{intthmnp}[thm]{Theorem}
\newenvironment{thmnp}{\begin{intthmnp}}{\qed \end{intthmnp}}

\theoremstyle{definition}
\newtheorem{defn}[thm]{Definition}
\newtheorem{warn}[thm]{Warning}

\theoremstyle{remark}
\newtheorem{rem}[thm]{Remark}

\def\thmref#1{Theorem~\ref{#1}}
\def\corref#1{Corollary~\ref{#1}}
\def\defref#1{Definition~\ref{#1}}
\def\eqref#1{Equation~\ref{#1}}
\def\ineqref#1{Inequality~\ref{#1}}
\def\figref#1{Figure~\ref{#1}}
\def\remref#1{Remark~\ref{#1}}
\def\secref#1{Section~\ref{#1}}

\def\co{\colon\thinspace}

\title{A note on pseudo-Anosov maps with small growth rate}
\date{{}}
\author{Peter Brinkmann}

\begin{document}
\maketitle

\begin{abstract}
We present an explicit sequence of pseudo-Anosov maps
$\phi_k: S_{2k}\rightarrow S_{2k}$ of surfaces of genus $2k$ whose growth
rates converge to one.
\end{abstract}

\section*{Introduction}

In this note, we present an explicit sequence $\phi_k$ of pseudo-Anosov
of surfaces of genus $2k$ whose growth rates converge
to one. This answers a question of Joan Birman,
who had previously asked whether such growth rates are bounded away
from one.  Norbert A'Campo, Mladen Bestvina, and Klaus Johannson
independently communicated this question to me. McMullen previously
obtained a similar result using quite different techniques \cite{mcm}.

The growth of the genus is not an artifact of our construction. For a
surface $S$ of fixed genus $g$, the growth rates of pseudo-Anosov maps
of $S$ are clearly bounded away from one, for they are
Perron-Frobenius eigenvalues of irreducible integral $m\times m$ matrices,
with $m\leq 6g-3$ \cite{hb2}. Finding the smallest possible growth rate
for each genus is an interesting problem that remains open.

One curious observation, due to Norbert A'Campo, is that our sequence
of pseudo-Anosov maps is a sequence of monodromies of w-slalom knots,
as defined in \cite{acampo1}.

In \secref{train}, we review the part of the theory of train tracks
\cite{hb1,hb2} that we use in this paper. \secref{exp} explains the
intuition that led to the result, and \secref{mainsec} contains the
statement and proof of the main results (\thmref{mainthm} and
\corref{maincor}).

The result of this paper grew out of massive computer experiments with
my software package XTrain \cite{pbexp,pbss} in the context of the REU
program at the University of Illinois at Urbana-Champaign.
I would like to thank Vamshidhar Kommineni for collecting much of the
experimental data that started this project. I am indebted to the
Department of Mathematics at UIUC for funding the computer experiments.
Finally, this paper would not have existed if Saul Schleimer had
not encouraged me to write it up.

\section{Train tracks}\label{train}

We present a brief review of train tracks as defined in \cite{hb1}.
Let $G$ be a finite graph without vertices of valence one or two, and
let $f\co G\rightarrow G$ be a homotopy equivalence of $G$ that maps
vertices to vertices. The map $f$ is said to be a {\em train track map}
if for every integer $n\geq 1$ and every edge $e$ of $G$, the restriction
of $f$ to the interior of $e$ is an immersion.

If $E_1,\cdots,E_m$ is the collection of edges of $G$, the
{\it transition matrix} of $f$ is the nonnegative
$m\times m$-matrix $M$ whose $ij$-th
entry is the number of times the $f$-image of $E_j$ crosses $E_i$, regardless
of orientation. $M$ is said to be {\it irreducible} if
for every tuple $1\leq i,j \leq m$, there exists some exponent $n>0$ such that
the $ij$-th entry of $M^n$ is nonzero.
If $M$ is irreducible, then it has a maximal real eigenvalue
$\lambda\geq 1$ (see \cite{seneta}).  We call $\lambda$ the
{\it growth rate} of $f$.

The following theorem from \cite{hb1} will be our main tool. Recall
that an outer automorphism $\omega$ of a free group $F$ is called
{\em reducible} if there are proper free factors $F_1,\ldots, F_r$
of $F$ such that $\omega$ permutes the conjugacy classes of the
$F_i$s and $F_1*\cdots*F_r$ is a free factor of $F$; $\omega$ is
{\em irreducible} if it is not reducible. Also, note that $\pi_1G$
is a finitely generated free group, and that a homotopy equivalence
$f\co G\rightarrow G$ induces an outer automorphism of $\pi_1G$.

\begin{thmnp}[{\cite[Theorem 4.1]{hb1}}]\label{maincrit}
Let $\omega$ be an outer automorphism of a finitely generated free
group $F$. Suppose that each positive power of $\omega$ is irreducible
and that there is a nontrivial word $s\in F$ such that $\omega$ 
preserves the conjugacy class of $s$ (up to inversion). Then
$\omega$ is geometrically realized by a pseudo-Anosov homeomorphism
$\phi\co S\rightarrow S$ of a surface with one puncture.
\end{thmnp}

\begin{rem}\label{maincrit2}
If $f\co G\rightarrow G$ is a train track map that induces an outer
automorphism $\omega$ as in \thmref{maincrit}, then
the transition matrix of $f$ is irreducible, and the growth rate of $f$ is the
same as the pseudo-Anosov growth rate of $\phi$.

Moreover, if $f\co G\rightarrow G$ is a train track map such that all
positive powers of its transition matrix $M$ are irreducible, then all
positive powers of the induced outer automorphism $\omega$ are
irreducible \cite{hb1}.
\end{rem}

\begin{rem}\label{remfol}
The proof of \corref{maincor} uses an explicit construction of invariant
foliations for pseudo-Anosov maps. This construction is straightforward
but too long to be reviewed in this note, so that we just point the reader
to \cite{hb2} for details.
\end{rem}

\section{Motivation}\label{exp}

\begin{warn} The discussion in this section is not supposed to present
any rigorous mathematical reasoning. Rather, the purpose of this section
is to explain the origin of the technical definitions and computations
of \secref{mainsec}.
\end{warn}

\begin{figure}
\renewcommand{\epsfsize}[2]{0.6\textwidth}
\psfrag{a0}{$a_0$}
\psfrag{b0}{$b_0$}
\psfrag{c0}{$c_0$}
\psfrag{d0}{$d_0$}
\psfrag{a1}{$a_1$}
\psfrag{b1}{$b_1$}
\psfrag{c1}{$c_1$}
\psfrag{d1}{$d_1$}
\centerline{\epsfbox{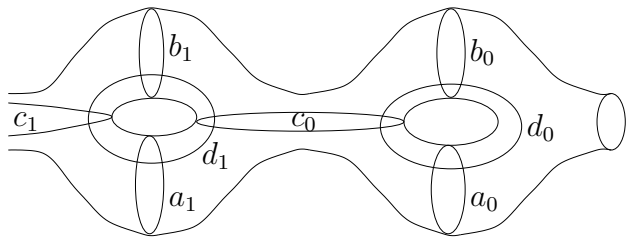}}
\caption{Generators of the mapping class group}\label{mcggen}
\end{figure}

One crucial tool in the development of the intuition behind \thmref{mainthm}
was XTrain \cite{pbexp,pbss}, a software package that implements algorithms
from \cite{hb1,hb2}, among others. In particular, the software allows users
to define homeomorphisms of surfaces with one puncture as a composition
of Dehn twists with respect to the curves shown in \figref{mcggen}.
When computing Dehn twists, we adopt the following convention: We equip the
surface with an outward pointing normal vector field. When twisting with
respect to a curve $c$, we turn {\it right} whenever we hit $c$. We denote
by $D_c$ the twist with respect to $c$.

The software represents a surface homeomorphism $\phi$ of a punctured surface
$S$ as a homotopy equivalence $f$ of a graph $G$ that is embedded in as well
as homotopy equivalent to $S$. There exists a loop $\sigma$ in $G$ that
corresponds to a short loop around the puncture of $S$. In particular,
$f$ preserves the free homotopy class of $\sigma$ (up to orientation).

\begin{figure}
\renewcommand{\epsfsize}[2]{0.5\textwidth}
\psfrag{x0}{$x_0$}
\psfrag{x1}{$x_1$}
\psfrag{x2}{$x_2$}
\psfrag{x2g+1}{$x_{2g}$}
\psfrag{vdots}{$\vdots$}
\center{\epsfbox{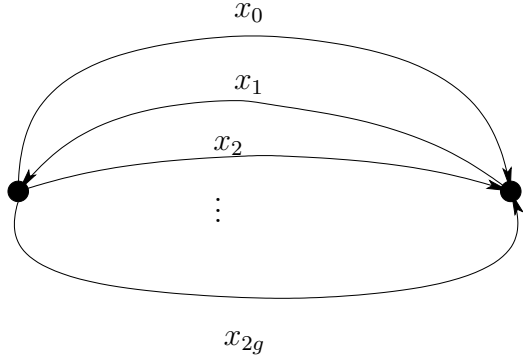}}
\caption{The graph $H_g$.}
\label{defh}
\end{figure}

The first ingredient is the observation that a homeomorphism of a surface of
genus $g$ given by
\begin{equation}\label{periodic}
\phi_g=D_{c_0}\cdots D_{c_{g-1}}D_{d_0}\cdots D_{d_{g-1}}
\end{equation}
can be represented by a train track map of a graph $H_g$ as in \figref{defh},
such that $x_0\mapsto x_1, x_1\mapsto x_2, \ldots, x_{2g}\mapsto x_0^{-1}$
with $\sigma_g=x_0x_1\cdots x_{2g}x_0^{-1}x_1^{-1}\cdots x_{2g}^{-1}$.
Note, in particular, that this map cyclically permutes the edges of $H_g$
(up to orientation).

The second ingredient comes from certain PV-automorphisms $\psi_n$
\cite{PV0} of a free group $F=\langle y_0,\ldots,y_n\rangle$ given by
$y_0\mapsto y_1, y_1\mapsto y_2,\ldots, y_n\mapsto y_0y_1$.
Mathematically, these automorphisms are very different from the maps
constructed in the previous paragraph (after all, PV automorphisms are
nongeometric and of exponential growth, whereas the maps of the previous
paragraph are geometric and periodic).

Superficially, though, these two classes of maps look strikingly similar.
Moreover, the growth rates of the maps $\psi_n$ converge to one. These two
observations prompted me to investigate maps that are built from blocks as
in \eqref{periodic}. Maps of surfaces of genus $2k$ of the form
\[
\phi_k=D_{c_0}\cdots D_{c_{k-1}}D_{d_0}\cdots D_{d_{k-1}}(D_{c_k}\cdots
D_{c_{2k-1}}D_{d_k}\cdots D_{d_{2k-1}})^{-1}
\]
turned out to be pseudo-Anosov with rather small growth rate. Computer
experiments suggested that the growth rates of these maps converges to
one, and the same experiments suggested that train tracks representing
these maps conform to a describable pattern, which gave rise to
\protect{\defref{maindef}} and \thmref{mainthm}. Notice how
\protect{\defref{maindef}} seems
reminiscent of both PV automorphisms as well as homeomorphisms as in
\eqref{periodic}.

\section{The sequence}\label{mainsec}

Motivated by the discussion of \secref{exp}, we now define a sequence of
surface homeomorphisms.

\begin{figure}
\renewcommand{\epsfsize}[2]{0.6\textwidth}
\psfrag{a}{$a$}
\psfrag{b}{$b$}
\psfrag{c}{$c$}
\psfrag{d}{$d$}
\psfrag{x0}{$x_0$}
\psfrag{x1}{$x_1$}
\psfrag{x2}{$x_2$}
\psfrag{x3}{$x_3$}
\psfrag{x2k-2}{$x_{2k-2}$}
\psfrag{x2k-1}{$x_{2k-1}$}
\psfrag{y0}{$y_0$}
\psfrag{y1}{$y_1$}
\psfrag{y2}{$y_2$}
\psfrag{y3}{$y_3$}
\psfrag{y2k-2}{$y_{2k-2}$}
\psfrag{y2k-1}{$y_{2k-1}$}
\psfrag{ddd}{$\cdots$}
\psfrag{vddd}{$\vdots$}
\center{\epsfbox{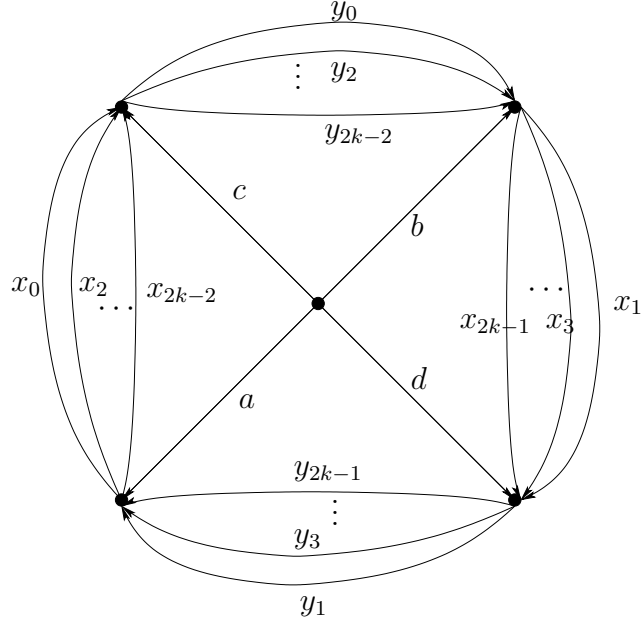}}
\caption{The graph $G_k$.}
\label{defg}
\end{figure}

\begin{defn}\label{maindef}
Let $k\geq 1$ be an integer, and choose the graph $G_k$ be as in \figref{defg}. 
We define a map $f_k\co G_k\rightarrow G_k$ by letting
\begin{align*}
a &\mapsto ax_0y_0\\
b &\mapsto by_0^{-1}x_0^{-1}\\
c &\mapsto d\\
d &\mapsto dy_1x_0\\
x_0 &\mapsto x_1\\
x_1 &\mapsto x_2\\
&\vdots\\
x_{2k-1} &\mapsto a^{-1}by_0^{-1}\\
y_0 &\mapsto y_1\\
y_1 &\mapsto y_2\\
&\vdots\\
y_{2k-1} &\mapsto c^{-1}b.\\
\end{align*}

Finally, let
\begin{align*}
\sigma_k=&x_0y_0x_1y_1\cdots x_{2k-1}y_{2k-1}a^{-1}by_0^{-1}c^{-1}d\\
&x_{2k-1}^{-1}b^{-1}cx_{2k-2}^{-1}y_{2k-1}^{-1}x_{2k-3}^{-1}y_{2k-2}^{-1}\cdots x_0^{-1}y_1^{-1}d^{-1}a.
\end{align*}
\end{defn}

We are now ready to state and prove the main result of this note.

\begin{thm}\label{mainthm}
The sequence of maps $f_k\co G_k\rightarrow G_k$ is a sequence of
homotopy equivalences induced by pseudo-Anosov maps
$\phi_k\co S_{2k}\rightarrow S_{2k}$ of surfaces of genus $2k$ with one
puncture. If $\lambda_k$ is the pseudo-Anosov growth rate of $\phi_k$,
then
\[
\lim\limits_{k\rightarrow \infty} \lambda_k = 1.
\]
\end{thm}

\begin{proof}
A number of tedious but straightforward checks yields the following
facts:

\begin{enumerate}
\item The maps $f_k$ are train track maps.
\item All positive powers of the transition matrix $M_k$ of $f_k$ are
irreducible.
\item The map $f_k$ preserves the free homotopy class of the loop
$\sigma_k$.
\end{enumerate}

Hence, by \thmref{maincrit} and \remref{maincrit2}, the outer
automorphism induced by $f_k$ is induced by a pseudo-Anosov map
$\phi_k\co S_k\rightarrow S_k$, and a quick computation of Euler
characteristics shows that the genus of $S_k$ is $2k$. Finally,
a simple induction shows that the characteristic polynomial of
the transition matrix $M_k$ is of the form
\[
\chi(\lambda)=
(\lambda-1)^2(\lambda^{4k+2}-\lambda^{4k+1}-4\lambda^{2k+1}-\lambda+1).
\]

Solving for the growth rate $\lambda_k$, we obtain
\begin{equation}\label{charpoly}
\lambda_k=1+\lambda_k^{4k+2}-\lambda_k^{4k+1}-4\lambda_k^{2k+1}.
\end{equation}

Note that the polynomial $\chi$ is palindromic (this is no surprise as
$f_k$ is induced by a surface homeomorphism), i.e.,
$\chi(\lambda)=\lambda^{4k+4}\chi(\frac{1}{\lambda})$.
Hence, \eqref{charpoly} also holds for $\lambda_k^{-1}$:
\begin{align}
\lambda_k^{-1}&=
	1+\lambda_k^{-(4k+2)}-\lambda_k^{-(4k+1)}-4\lambda_k^{-(2k+1)}\nonumber\\
&\geq 1-\lambda_k^{-(4k+1)}-4\lambda_k^{-(2k+1)}.\label{mainineq}
\end{align}

Recall that $\lambda_k^{-1}<1$. Let $0<u<1$ be some real number. We have
$\lim_{k\rightarrow \infty} 1-u^{4k+1}-4u^{2k+1}=1$, which implies that
$u$ only satisfies \ineqref{mainineq} for finitely many values of $k$.
Hence, for any such $u$, the set $\{\lambda_k | \lambda_k^{-1}<u\}$ is
finite. This immediately implies that
$\lim_{k\rightarrow \infty} \lambda_k^{-1}=1$, hence
\[
\lim_{k\rightarrow \infty} \lambda_k=1.
\]
\end{proof}

\begin{cor}\label{maincor}
The maps $\phi_k\co S_{2k} \rightarrow S_{2k}$ from \thmref{mainthm}
can be extended to pseudo-Anosov maps of closed surfaces. The growth
rates of the extended maps are the same as those of the original maps.
\end{cor}

\begin{proof}
A lengthy but straightforward computation of invariant foliations
(see \remref{remfol}) yields that the four outer vertices of the graph
in \figref{defg} give rise to singularities of index $\frac{1}{2}-k$,
while the central vertex does not give rise to any singularity. Hence,
the sum of the indices of all singularities coming from vertices of
the graph is $2-4k$, which is the Euler characteristic of a closed
surface of genus $2k$.

Hence, the foliations have no singularity at the puncture, which
implies that the extension of $\phi_k$ to the closed surface obtained
by filling in the puncture is pseudo-Anosov, with the same growth
rate as $\phi_k$.
\end{proof}

\bibliographystyle{alpha}
\bibliography{all}
\par

{\noindent \sc Department of Mathematics\\
University of Illinois at Urbana-Champaign\\
273 Altgeld Hall\\
1409 W.\ Green St.\\
Urbana, IL 61801, USA\\}
\par
{\noindent \it E-mail:} brinkman@math.uiuc.edu

\end{document}